\documentclass[12pt]{amsart}
\usepackage{latexsym}
\usepackage{a4}
\usepackage{amssymb}
\usepackage{amsbsy}

\newtheorem{thm}{Theorem}[section] 
 
\newtheorem{lem}[thm]{Lemma} 
 
\newtheorem{cor}[thm]{Corollary} 
\theoremstyle{definition}
 
\newtheorem{de}[thm]{Definition}

\numberwithin{equation}{section} 
\newcommand{\ab}[1]{{\mathbf{#1}}}

\newcommand{\ul}[1]{\{1,\dots,#1\}}

\newcommand{\N}{\Bbb{ N}}

\newcommand{\setsuchthat}{\,\, \pmb{|} \,\,}

\newcommand{\vb}[1]{\mathbf{#1}} 
\newcommand{\cb}[1]{#1}

\newcommand{\Pol}{\mathrm{Pol}\,}

\newcommand{\A}{A}

\newcommand{\DCC}{{\mathrm{(DCC)}}}

\newcommand{\ltlex}{<_{\mathrm{lex}}}
\newcommand{\lelex}{\le_{\mathrm{lex}}}
\newcommand{\lea}{\le_{\mathrm{E}}}

\newcommand{\FO}{\textsf{firstOcc}\,}
\newcommand{\Tab}{T_{\vb{a},\vb{b},h}}
\newcommand{\ran}{\textrm{range}}
\newcommand{\Symbols}{\textsf{Symbols}\,}
\newcommand{\Start}{\textsf{Start}\,}
\newcommand{\Last}{\textsf{Last}\,}

\newcommand{\algop}[2]{\langle {#1}, {#2} \rangle}

\newcommand{\PhiX}{\varphi}
\newcommand{\Sig}{\mathrm{Sig}}
\title{On the number of finite algebraic structures}
\author{Erhard Aichinger} 
\address{Erhard Aichinger,
Institut f\"ur Algebra,
Johannes Kepler Universit\"at Linz,
4040 Linz,
Austria}
\email{\tt erhard@algebra.uni-linz.ac.at}
\author{Peter Mayr} 
\address{Peter Mayr,
Centro de \'Algebra da Universidade de Lisboa (CAUL),
 1649-003 Lisboa, Portugal \&
Institut f\"ur Algebra,
Johannes Kepler Universit\"at Linz,
4040 Linz,
Austria
}
\email{stein@cii.fc.ul.pt}
\author{Ralph McKenzie}
\address{Ralph McKenzie,
Department of Mathematics,
Vanderbilt University,
Nashville, Tennessee, U.S.}
\email{ralph.n.mckenzie@vanderbilt.edu}
\date{November 30, 2010}

\begin{document}
\bibliographystyle{amsalpha}

\maketitle
\begin{abstract}
 We prove that every clone of operations on a finite set $A$, if it
 contains a Malcev operation, is finitely related -- i.e., identical with the
 clone of all operations respecting $R$ for some finitary 
  relation $R$ over $A$. It follows that for a fixed finite set
  $A$, the set of all such Malcev clones is countable.
 This completes the solution of a problem that was first formulated in 1980,
 or earlier: how many Malcev clones can finite sets support? More generally,
 we prove that every finite algebra with few subpowers has a finitely related
 clone of term operations.
  Hence modulo term equivalence and a renaming of the elements, 
  there are only countably
  many finite algebras with few subpowers, and thus only
  countably many finite algebras with a Malcev term.
\end{abstract}
   
\section{Introduction}

 An algebraic structure (or \emph{algebra}, for short) is usually
 represented as a nonvoid set together with a 
 set of finitary operations on it.
 In the present paper,
 we contribute to the following question: how many essentially different finite
 algebraic structures exist? 
 Clearly, on
 a finite set of size at least two, there are countably many 
 finitary operations, and hence there are continuum
 many ways to choose a set of basic operations.
However, many of these algebras
are equivalent in the sense that the same functions can be 
composed from their basic operations; these compositions
are called the \emph{term functions} of the algebra.
 Two algebras
are \emph{term equivalent} if they have the same set of 
term functions. The Boolean algebra
 $\algop{B}{\wedge,\vee,\neg}$ and its 
counterpart, the Boolean ring $\algop{B}{+,\cdot, 1}$, are examples
of term equivalent algebras.
Many structural properties of an algebra,
like its subalgebras, congruence relations, automorphisms, etc.,
 depend on its term functions rather than on
the particular choice of basic operations. 
 Hence we are motivated to classify algebras modulo term equivalence.
 In 1941 E.~Post~\cite{Po:TTVI} published that there are only countably many term inequivalent
 algebras of size two (modulo renaming of the elements),
 and he described them all explicitly.
 In 1959 J.~Janov and A.~Mu{\v{c}}nik~\cite{JM:EKCC}
 showed that even modulo term
equivalence,
the number of algebras on a finite set with 
at least three elements is uncountable.

 Many classical algebraic structures
 have the property that their congruence
 relations commute with respect to the relation product.
 A.~Malcev~\cite{Ma:OTGT} has characterized varieties of algebras with this property
 (a variety is a class of algebras of the same type that is defined by equations);
 a consequence of his result is that an algebra generates such a
 congruence-permutable variety if and only if it has a ternary (Malcev) term
 operation $m$ satisfying $m(x,y,y)= m(y,y,x) = x$ for all $x,y$.  
 These algebras include
all finite algebras that have a quasigroup operation
among their binary term functions, and hence, e.g.,
all finite groups, rings, modules, loops, and planar ternary rings.
 It has long been open how many of the $2^{\aleph_0}$ finite term inequivalent algebras
 on a set of size at least three have a Malcev term
 (see e.g.~\cite[Problem 5.19]{KP:PARI}). We will prove that this number is
 at most countably infinite.
 In particular, Theorem~\ref{thm:countable} yields that for every finite algebra
 $\ab{A}$ with a Malcev term there is an
 $n \in \N$ and a single subalgebra $R$ of $\ab{A}^n$ such that $\ab{A}$ is determined by
 $R$ up to term-equivalence.

 Recently a combinatorial characterization of finite algebras with a Malcev term
 has been found. As a consequence
of \cite{BI:VWFS}, a finite algebra
$\ab{A}$ has a Malcev term if and only if there is 
a positive real $c$ such that every independent subset of $\ab{A}^n$ has at most
$cn$ elements (Here a subset $X$ is independent if no proper subset of $X$ generates
 the same subalgebra of $\ab{A}^n$ as $X$).
 This condition immediately yields that $\ab{A}^n$ has at most $|A|^{cn^2}$ 
subalgebras.
 In general, a finite algebra $\ab{A}$ for which there exist a polynomial $p$ such
 that $\ab{A}^n$ has at most $2^{ p(n)}$ subalgebras is
 said to have \emph{few subpowers}
 (Note that the number of subalgebras of $\ab{A}^n$ is certainly bounded by $2^{|A|^n}$.
 The adjective `few' refers to the fact that the number of subalgebras does not grow doubly
 exponential in $n$). 
 In~\cite{BI:VWFS} algebras with few subpowers
are characterized by the existence of an \emph{edge operation}
(see Section~\ref{sec:ac}) 
 among their term functions.
The class of algebras with an edge term
is a vast extension of the class
of algebras with a Malcev term. It also comprises, e.g., 
all lattices and algebras with lattice operations,
and is properly contained in the class of algebras that generate congruence modular
varieties. 
 Theorem~\ref{thm:countable} yields that every finite algebra with few subpowers
 is finitely related (see Section~\ref{sec:ac}).
 This means that every such algebra -- even if it has an infinite set of basic operations --
 has a finite description up to term equivalence.
Hence on a finite set $A$, 
modulo term equivalence, the number of
algebras with few subpowers is at most countably infinite (Corollary~\ref{cor:countable}).

 Algebras with few subpowers
 recently appeared in connection with the constraint satisfaction problem
 (CSP) in computer science. By~\cite{IMMVW:TALA} CSPs that afford an edge term
 can be solved by a polynomial-time algorithm.
 It is expected that more generally, CSPs admissible over finite algebras in
 congruence-modular varieties are solvable in polynomial time as well.
 This would follow from a partial converse of our result which has been
 conjectured by M.\ Valeriote. The conjecture is that a finite algebra in a
 congruence-modular variety, if it is finitely related,
 must have few subpowers. 
 A special case of this, which had earlier been conjectured by
 L.\ Z{\'a}dori, has been established recently by L. Barto~\cite{Ba:CDINU}
 (see also P. Markovi\'{c} and R. McKenzie~\cite{MM:FSCD}): A finite algebra in a
 congruence-distributive variety is finitely related if and only if it has a
 near-unanimity operation.

\section{Algebras and Clones} \label{sec:ac}

    We will express our results using the terminology of
    universal algebra \cite{BS:ACIU,MMT:ALVV} and 
    clone theory \cite{PK:FUR,Sz:CIUA}.
    Following \cite{HM:TSOF},
    we understand an algebra $\ab{A} := \langle A,F \rangle$ as a set $A$ together with 
    a set of finitary operations $F$ on $A$.
    For a non-void set $A$, by a \emph{clone} on $A$
    we shall mean any set of finitary operations on $A$ (of positive arity)
    that
    is closed under compositions and contains the projection
    operations $e_i^n (x_1,\ldots,x_n) = x_i$ for all positive integers $n$
 and for all $i \in \{1,\ldots, n\}$. 
    The set of term operations of an algebra $\ab{A}$
    is a clone, and every clone
    on $A$ takes this form.
  
 For $k\geq 2$ a function $t: A^{k+1}\rightarrow A$ is a
 $k$-\emph{edge operation} if for all $x,y\in A$ we have
\[ t(y,y,x,\ldots,x) = t(y,x,y,x,\ldots,x) = x \]
 and for all $i \in \{4,\ldots, k+1\}$ and for all $x,y \in A$, we have
\[ t(x, \ldots,x, y, x, \ldots, x) = x, \text{ with } y \text{ in position } i. \]
 A ternary operation $t$ is a $2$-edge operation if and
 only if $m(x,y,z) := t(y,x,z)$ is a Malcev operation.
 For $k > 2$ a $k$-ary \emph{near unanimity} operation $f$ is a function
 such that $t(x_1,\ldots,x_{k+1}) := f(x_2,\ldots,x_{k+1})$ is a $k$-edge operation.
 Thus the class of clones with edge operations contains all clones with Malcev
 or near unanimity operations.
 We also note that
 an algebra has an edge term if and only if it has a \emph{parallelogram term} as defined
 in~\cite{KS:CAWP}.

 A clone $C$ on $A$ is \emph{finitely related} if there exist subalgebras
 $R_1,\ldots,R_k$ of finitary powers of $\algop{A}{C}$ such that every function
 on $A$ that preserves every $R_i$ for $i\in\{1,\ldots,k\}$ is in $C$.
 We call an algebra \emph{finitely related} if its clone of term functions is
 finitely related.
 Clones containing a near-unanimity operation are finitely related by the Baker-Pixley
 Theorem \cite{BP:PIAT}.
 In~\cite{Ai:CMCO} the first author shows that, on a finite set, every clone that
 contains a Malcev operation and all constant functions, is finitely related. 
 Special cases of the result in~\cite{Ai:CMCO} were given, for example, by 
 P.~Idziak~\cite{Id:CCMO},
 A.~Bulatov~\cite{Bu:OTNO},
 K.~Kearnes and \'A.~Szendrei~\cite{KS:COFG},
 the second author~\cite{Ma:PCOS, Ma:MASC}, N.~Mudrinski and the first author
 \cite{AM:PCOM}.
 In this paper we prove the common generalization that on a finite set
 every clone with edge operation is finitely related (Theorem~\ref{thm:relations}).

  The conjecture that on a finite set the number of clones with Malcev operation
  is countable dates back to the mid 1980's or earlier. 
  The two tools which we use to prove this conjecture were first
   combined to good effect in \cite{Ai:CMCO}. They are, first,
  a combinatorial theorem due to G.\ Higman \cite{Hi:OBDI}, which occurs
  here in a generalized form as Lemma~\ref{lem:leawpo}; and second, the result that
 for an algebra $\A$ with $k$-edge term every subalgebra of a finite power of $\A$
 has a small  generating set that takes a specific form (Lemma~\ref{lem:rep}).  
  The second result also lies at the core of the proof in \cite{IMMVW:TALA}
  that every constraint satisfaction problem whose template relations
  are admissible over an algebra with few subpowers, is
  tractable -- i.e, admits a polynomial time algorithm for its
  solution.

\section{Preliminaries from order theory}   \label{sec:order}

We will first give a short survey of those results from
order theory that we will need in the sequel.
The partially ordered set $\algop{X}{\le}$ is
\emph{well partially ordered} if it satisfies the descending chain
condition $\DCC$
and has no infinite antichains.
The following facts about well partial orders can be found in \cite{La:WASO}
(cf. \cite{Na:OFT}).
A sequence of elements $\langle x_k \setsuchthat k \in \N \rangle$
is \emph{good} if there are $i,j \in \N$ with $i < j$ and
$x_i \le x_j$; a sequence is \emph{bad} if it is not good.
Using Ramsey's Theorem, one can prove
that $\algop{X}{\le}$ is well partially ordered if and only if
every sequence in $X$ is good.
If $\algop{X}{\le}$ satisfies the $\DCC$, but is not well partially ordered,
then there exists a bad sequence 
$\langle x_k \setsuchthat k \in \N \rangle$ with the property
that for all $i \in \N$ and for all $y_i \in X$ with $y_i < x_i$,
every sequence starting with $(x_1, \ldots, x_{i-1}, y_i)$ is good.
Such a sequence is called a \emph{minimal bad sequence}.
For an ordered set $\algop{X}{\le}$, a subset $Y$ of $X$
is \emph{upward closed} if for all $y \in Y$ and $x \in X$ with
$y \le x$, we have $x \in Y$.

For $A = \{1, 2, \ldots, t \}$,  
we will use the lexicographic ordering
on $A^n$. 
For $\vb{a} = (a_1,\ldots, a_n)$ and
$\vb{b} = (b_1, \ldots, b_n)$, we say
$\vb{a} \lelex \vb{b}$ if
\begin{multline*}
     (\exists i \in \{1,\ldots,n\} :
       a_1 = b_1 \wedge \ldots \wedge
       a_{i-1} = b_{i-1} \wedge a_i < b_i) \text{ or }
      \\ (a_1, \ldots,a_n) = (b_1,\ldots, b_n).
\end{multline*}

For every finite set $A$,
we let $A^+$ be the set
$\bigcup \{A^n \setsuchthat n \in \N\}$.
We will now introduce an order relation on $A^+$.
For $\vb{a} = (a_1,\ldots, a_n) \in A^+$ and $b \in A$,
we define the \emph{index of the first occurrence of $b$ in $\vb{a}$},
$\FO (\vb{a}, b)$, by
$\FO (\vb{a}, b) := 0$ if $b \not\in \{a_1,\ldots, a_n \}$,
and $\FO (\vb{a}, b) := \min \{ i \in \{1,\ldots, n\} \setsuchthat
a_i = b\}$ otherwise.
\begin{de} Let $A$ be a finite set, and
    let $\vb{a} = (a_1,\ldots, a_m)$ and $\vb{b} = (b_1,\ldots, b_n)$ be
 elements of $A^+$. We say
    $\vb{a} \lea \vb{b}$
 (read: $\vb{a}$ embeds into $\vb{b}$)
 if there is an injective and
    increasing function $h :\{1,\ldots, m\} \to \{1,\ldots,n\}$
    such that
    \begin{enumerate}
         \item for all $i \in \{1,\ldots, m\}$ : 
                       $a_i = b_{h(i)}$,
         \item $\{a_1,\ldots, a_m\} = \{b_1,\ldots, b_n\}$,
         \item for all $c \in \{a_1,\ldots, a_m\}$:
                $h (\FO (\vb{a}, c)) = \FO (\vb{b}, c)$.
    \end{enumerate}
We will call such an $h$ a function \emph{witnessing $\vb{a} \lea \vb{b}$}.
\end{de}

 Less formally, we have $\vb{a} \lea \vb{b}$ for words $\vb{a},\vb{b}$ over the
 alphabet $A$ if and only if $\vb{b}$ can be obtained from $\vb{a}$ by inserting
 additional letters anywhere after their first occurrence in $\vb{a}$.  
 We will use the following fact about this ordering, which generalizes
 Higman's Theorem 4.4 in \cite{Hi:OBDI}.

\begin{lem} \label{lem:leawpo}
  Let $A$ be a finite set.
  Then $\algop{A^+}{\lea}$ is well partially ordered.
\end{lem}
\emph{Proof:}
 It is easy to see that $\lea$ is a partial order relation
 and that $\algop{A^+}{\lea}$ satisfies the $\DCC$.
 It remains to show that for every sequence
 $\langle \vb{x}^{(k)} \setsuchthat k \in \N \rangle$ in $A^+$, there exist 
 $i,j\in\N$ such that $i<j$ and $\vb{x}^{(i)}\lea \vb{x}^{(j)}$.
 We will prove this by induction on $|A|$.
 For $|A| = 1$, the claim is obvious.
 Assume $|A| > 1$ and that $\algop{B^+}{\lea}$ is well partially ordered for
 every proper subset $B$ of $A$. 

 Seeking a contradiction we suppose we have a minimal bad sequence
 $\langle \vb{x}^{(k)} \setsuchthat k \in \N \rangle$ in $A^+$.
    For each $\vb{x} = (x_1,\ldots,x_n) \in A^+$, let
    $\Symbols (\vb{x}) := \{ x_1, \ldots, x_n\}$ be the
    set of all elements of $A$ that occur in the word $\vb{x}$,
    let $\Last (\vb{x}) := x_n$ denote the last letter of $\vb{x}$,
    and, if $n \ge 2$, let
    $\Start (\vb{x}) := (x_1, \ldots, x_{n-1})$.
 Since $A$ is finite,
 we have $a \in A$ and an infinite $T \subseteq\N$ such that for all
 $i\in T$, $\Last (\vb{x}^{(i)}) = a$ and the length of $\vb{x}^{(i)}$
 is at least two.
 
 Let us first consider the case that there exist an infinite $S \subseteq T$
 such that $\Symbols (\Start (\vb{x}^{(i)})) \subseteq A \setminus \{a\}$ for
 all $i \in S$.
 By the induction hypothesis, $\lea$ is a well partial order on
 $(A \setminus \{a\})^+$. 
 Hence there are $i, j \in S$ with $i < j$ such that
 $\Start (\vb{x}^{(i)}) \lea \Start (\vb{x}^{(j)})$.
 Since $a$ does not occur in $\Start (\vb{x}^{(i)})$ nor in 
 $\Start (\vb{x}^{(j)})$, and since
 $\Last (\vb{x}^{(i)}) = \Last (\vb{x}^{(j)}) = a$, we have
 $\vb{x}^{(i)} \lea \vb{x}^{(j)}$, contradicting the fact that
 $\langle \vb{x}^{(k)} \setsuchthat k \in \N \rangle$ is a bad sequence.

 Thus we may assume that there exist an infinite subset
 $S := \{ s_1,s_2,\ldots \}$ of $T$ (with $s_i < s_j$ whenever $i<j$) such that
 $\Symbols (\Start (\vb{x}^{(s)})) = A$ for all $s\in S$.
 Now consider the sequence
\[ 
     \langle \vb{y}^{(k)} \setsuchthat k \in \N \rangle := 
      \langle \vb{x}^{(1)},\vb{x}^{(2)},\ldots,\vb{x}^{(s_1-1)},
              \Start(\vb{x}^{(s_1)}),\Start(\vb{x}^{(s_2)}),\ldots \rangle. 
 \] 
 We show that $\langle \vb{y}^{(k)} \setsuchthat k \in \N \rangle$ is bad
 by distinguishing three cases: 
 If $i < j < s_1$, then clearly $\vb{x}^{(i)}\not\lea\vb{x}^{(j)}$. 
 If $i < s_1$ and $j \geq 1$, then $\vb{x}^{(i)}\lea\Start(\vb{x}^{(s_j)})$
 yields $\vb{x}^{(i)}\lea\vb{x}^{(s_j)}$, contradicting
 the fact that $\langle \vb{x}^{(k)} \setsuchthat k \in \N \rangle$ is bad.
 If $i < j$, then $\Start(\vb{x}^{(s_i)})\lea\Start(\vb{x}^{(s_j)})$
 implies $\vb{x}^{(s_i)}\lea\vb{x}^{(s_j)}$ because 
 $\Last (\vb{x}^{(s_i)}) = \Last (\vb{x}^{(s_j)}) = a$ and $a$ already occurs
 both in $\Start(\vb{x}^{(s_i)})$ and in $\Start(\vb{x}^{(s_j)})$.
 This again contradicts the badness of
$\langle \vb{x}^{(k)} \setsuchthat k \in \N \rangle$.
 Hence $\langle \vb{y}^{(k)} \setsuchthat k \in \N \rangle$ is bad.
 However, since $\vb{y}^{(s_1)} = \Start(\vb{x}^{(s_1)}) <_E \vb{x}^{(s_1)}$,
 this contradicts the choice of
 $\langle \vb{x}^{(k)} \setsuchthat k \in \N \rangle$ as a minimal bad sequence.
 Hence $\algop{A^+}{\lea}$ is well partially ordered.
\qed

 For $\vb{a},\vb{b}\in A^+$ with $\vb{a}\lea\vb{b}$ we observe a correspondence
 between the elements that are lexicographically smaller than $\vb{a}$ and certain
 elements that are lexicographically smaller than $\vb{b}$. But before that
 we need to introduce some notation.

\begin{de} Let $A$ be a finite set, 
    let $\vb{a} = (a_1,\ldots,a_m) \in A^m$, 
    $\vb{b} = (b_1, \ldots, b_n) \in A^n$
    be such that $\vb{a} \lea \vb{b}$, and let $h$ be a function
    from $\{1,\ldots, m\} \to \{1,\ldots, n\}$ witnessing
    $\vb{a} \lea \vb{b}$.
    We 
 define a function $\Tab : A^m \to A^n$.
    Let $\vb{x} = (x_1, \ldots, x_m) \in A^m$.
    If $j \in \ran (h)$, then
    the $j$-th entry of $\Tab (\vb{x})$, abbreviated
    by $\Tab (\vb{x})\,(j)$, is defined by
    \[
       \Tab (\vb{x}) \, (j) := x_i,
    \]
     where $i \in \{1,\ldots, m\}$ is such that 
    $h(i) = j$.
    If $j \not\in \ran (h)$,
   then
   \[
      \Tab (\vb{x}) \, (j) := x_i,
   \]
    where $i := \FO (\vb{a}, b_j)$.
  \end{de}
      
 \begin{lem} \label{lem:co}
   Let $t \in \N$, let $A = \{ 1, 2, \ldots, t \}$,
   and let
   $\ab{a} \in A^m$, $\ab{b} \in A^n$ with
 $h: \{1,\ldots, m\} \to \{1,\ldots, n\}$ witnessing
    $\vb{a} \lea \vb{b}$.
   Let $\vb{c} \in A^m$ be such that
   $\vb{c} \ltlex \vb{a}$.
   Then we have
   \begin{enumerate} 
      \item \label{it:c1} $\Tab (\vb{a}) = \vb{b}$,
      \item \label{it:c2} $\Tab (\vb{c}) \ltlex \vb{b}$. 
   \end{enumerate}
\end{lem}
\emph{Proof:}
    \eqref{it:c1} follows immediately from the definition
    of $\Tab$. For proving~\eqref{it:c2},
    let $k$ be the index of the first place in which
    $\vb{c}$ differs from $\vb{a}$.
    Hence $\vb{c} = (a_1,\ldots, a_{k-1}, c_k, c_{k+1},\ldots)$,
    $\vb{a} = (a_1,\ldots, a_{k-1}, a_k, a_{k+1},\ldots)$, and
    $c_k < a_k$.
    
     We first show that for all $j < h(k)$, we have
    $\Tab (\vb{c}) (j) = \Tab (\vb{a}) (j)$.
     If $j$ is in the range of $h$, there is an $i$ with   
     $h(i) = j$, and we have
     $\Tab (\vb{c}) (j) = c_i$ and
     $\Tab (\vb{a}) (j) = a_i$.
     Since $h(i) < h(k)$, we have
     $i < k$. Thus $c_i = a_i$, since $k$ is the first index
     at which $\vb{c}$ and $\vb{a}$ differ.
     We now consider the case that
     $j$ is not in the range of $h$.
 Since $\{ b_1,\ldots, b_n\} = \{a_1,\ldots, a_m\}$, we have that
 $i := \FO (\vb{a}, b_j)$ satisfies $i > 0$. By the definition of $\lea$
 we have $h(i) = \FO (\vb{b}, b_j)$ and therefore $h(i) \le j$.
 Hence $h(i) < h(k)$ and $i < k$. Thus $c_{i} = a_{i}$.
 Since $\Tab (x_1, \ldots, x_m) (j) := x_{i}$ for all $\vb{x} \in A^m$, 
 we finally obtain $\Tab (\vb{c}) (j) = \Tab (\vb{a}) (j)$.
    
     Since $\Tab (\vb{a}) (h(k)) = a_k$ and
     $\Tab (\vb{c}) (h(k)) = c_k$, we have
     $\Tab (\vb{c}) \ltlex \Tab (\vb{a})$.
\qed

\section{Algebras with edge term} \label{sec:edge}

 Let $A$ be a set, and let $m\in\N$. For $\vb{a} = (a_1,\ldots,a_m) \in A^m$
 and $T\subseteq\{1,\ldots,m\}$, 
 we denote the projection to the tuple of entries that are indexed by $T$ as
\[ \pi_T(\vb{a}) := \langle a_i \setsuchthat i \in T \rangle. \]
 For $F \subseteq A^m$ and $i\in\{1,\ldots,m\}$, define
\[ \varphi_i(F) := \{ (a_i,b_i) \in A^2 \setsuchthat \vb{a},\vb{b}\in F \text{ and } \pi_{\{1,\ldots,i-1\}}(\vb{a}) = \pi_{\{1,\ldots,i-1\}}(\vb{b}) \}. \]
 By~\cite[Lemma 3.1]{Ai:CMCO} a subuniverse $G$ of a Malcev algebra $\ab{A}^m$
 is generated by every subset $F$ of $G$ with $\varphi_i(F) = \varphi_i(G)$
 for all $i\in\{1,\ldots,m\}$.

 In~\cite{BI:VWFS} these relations $\varphi_i$ and projections $\pi_T$ occur
 in the description of small generating sets for the subuniverses of $\ab{A}^m$
 for a finite algebra $\ab{A}$ with edge term operation. These generating sets
 were then used to obtain a bound on the number of subuniverses of $\ab{A}^m$.
 We reformulate the representation result~\cite[Corollary 3.9]{BI:VWFS} for
 our purposes.

\begin{lem}  \label{lem:rep}
 Let $k,m$ be positive integers with $k>1$,
 let $\ab{A}$ be a finite algebra with $k$-edge term operation $t$, and
 let $F,G$ be subuniverses of $\ab{A}^m$ with $F \subseteq G$.
 Assume $\pi_T(F) = \pi_T(G)$ for all $T\subseteq\{1,\ldots,m\}$ with $|T| < k$,
 and $\varphi_i(G) \subseteq \varphi_i(F)$ for all $i\in\{1,\ldots,m\}$.
 Then $F=G$.
\end{lem}
 
\emph{Proof:}
 We only have to check that $F$ is what is called a \emph{representation} of
 $G$ in~\cite[Def.~3.2]{BI:VWFS}. For that we let $d$ be the binary term
 function on $\ab{A}$ that is defined from $t$ in Lemma 2.13 of~\cite{BI:VWFS}.
 We also need the notion of a
 \emph{signature} $\Sig_R$ of a subset $R$ of $A^m$,
\[ \Sig_R := \{ (i,u,v) \in \{1,\ldots,m\}\times A^2 \setsuchthat
 (u,v)\in\varphi_i(R) \text{ and } d(u,v) = v \}. \]

 From $F \subseteq G$, it is immediate that
 $\varphi_i(F) \subseteq \varphi_i(G)$.
 Consequently $\varphi_i(F) = \varphi_i(G)$ for all $i\in\{1,\ldots,m\}$.
 In particular $\Sig_F = \Sig_G$.
 Thus $F$ is a representation of $G$. Since $F,G$ are subuniverses of
 $\ab{A}^m$,
 Corollary 3.9 of~\cite{BI:VWFS} yields $F = G$.
\qed

 The previous result has also been known in two special cases:
 For $\ab{A}$ with a $k$-ary near unanimity term it follows from the 
 Baker-Pixley Theorem~\cite{BP:PIAT}. 
 For $\ab{A}$ with a Malcev term, it occurs as Lemma 3.1 in~\cite{Ai:CMCO},
 and it is the central fact underlying Dalmau's polynomial-time algorithm
 for solving CSPs 
 which admit a Malcev polymorphism~\cite{BD:ASAF}.

\section{Encoding clones}

Let $\cb{C}$ be a clone on the $t$-element set 
$A = \{1, 2, \ldots, t\}$, and let $n \in \N$.
 Let $\cb{C}^{[n]}$ denote the set of $n$-ary functions in $\cb{C}$.
 As in~\cite{Ai:CMCO}, for $\vb{a} \in A^n$, we define a binary 
relation $\PhiX (\cb{C}, \vb{a})$ on $A$ by
\[
    \PhiX (\cb{C}, \vb{a}) :=
       \{ (f(\vb{a}), g(\vb{a})) \setsuchthat
          f,g \in \cb{C}^{[n]},
          \forall \vb{c} \in A^n: 
            \vb{c} \ltlex \vb{a} 
             \Rightarrow f (\vb{c}) = g(\vb{c}) \}.
\]
 Intuitively, if $\PhiX (\cb{C}, \vb{a})$ is small, then the functions in $\cb{C}$
 are strongly restricted by their images on $\vb{c}$ for $\vb{c} \ltlex \vb{a}$.
 We also encode these relations in another way.

 For $(c,d)\in A^2$, we define a subset $\lambda (\cb{C}, (c,d))$ of $A^+$ by
   \[
      \lambda(\cb{C}, (c,d)) :=
        \{ \vb{a} \in A^+ \setsuchthat
           (c,d)\not\in \PhiX (\cb{C}, \vb{a}) \}.
   \]
 From the order theoretic observations in Section~\ref{sec:order} we obtain the following 
 lemmas.

\begin{lem} \label{lem:embeddingPhi}
    Let $t,m,n \in \N$, let $\cb{C}$ be a 
    clone on the $t$-element set $A =\{1,2,\ldots, t\}$,
    and let $\vb{a} \in A^m$, $\vb{b} \in A^n$ such that
    $\vb{a} \lea \vb{b}$.
    Then $\PhiX (\cb{C}, \vb{b}) \subseteq \PhiX (\cb{C}, \vb{a})$.
\end{lem}
\emph{Proof:}
   Let $(x,y) \in \PhiX (\cb{C}, \vb{b})$. Then there are
   $f,g \in \cb{C}^{[n]}$ such that $x = f(\vb{b})$,
   $y = g(\vb{b})$, and $f(\vb{c}) = g (\vb{c})$ for all
   $\vb{c} \in A^n$ with $\vb{c} \ltlex \vb{b}$.
   Let $h$ be a function from $\{1,\ldots, m\}$ to $\{1,\ldots,n\}$
   witnessing $\vb{a} \lea \vb{b}$.
   Now we define functions $f_1$ and $g_1$ from $A^m
$ to $A$ by
 \[
     \begin{array}{rcl}
        f_1 (\vb{x}) & := & f (\Tab (\vb{x})) \\
        g_1 (\vb{x}) & := & g (\Tab (\vb{x}))
     \end{array}
 \]
 for $\vb{x} \in A^m$.
 By the definition of $\Tab$, we see that for each $j \in \{1,\ldots, n\}$,
 the mapping that maps $\vb{x}$ to the $j$-th component of
 $\Tab(\vb{x})$ is a projection operation.
  Hence $f_1$ and $g_1$ lie in the clone $\cb{C}$.
 
We will now show that $(f_1 (\vb{a}), g_1 (\vb{a}))$ is an element of 
$\PhiX (\cb{C}, \vb{a})$. To this end,
 let $\vb{c} \in A^m$ be such that $\vb{c} \ltlex \vb{a}$.
 Then Lemma~\ref{lem:co} yields $\Tab (\vb{c}) \ltlex \vb{b}$.
 Hence we have
 $f_1 (\vb{c})= f (\Tab (\vb{c})) = g (\Tab (\vb{c})) = g_1 (\vb{c})$.
From this we obtain $(f_1 (\vb{a}), g_1 (\vb{a})) \in \PhiX (\cb{C}, \vb{a})$.
Since $(f_1 (\vb{a}), g_1 (\vb{a})) = (f (\vb{b}), g (\vb{b}))
        = (x,y)$ by Lemma~\ref{lem:co},
 we obtain $(x,y) \in \PhiX (\cb{C}, \vb{a})$.
\qed
                                  
\begin{lem} \label{lem:upwardclosed}
  Let $\cb{C}$ be a clone on a finite set $A$, and let $(c,d)\in A^2$.
  Then $\lambda (\cb{C}, (c,d))$ is an upward closed subset
  of $\algop{A^+ }{\lea}$.
\end{lem}
\emph{Proof:}
    Let $\vb{a} \in \lambda (\cb{C}, (c,d))$, and let
    $\vb{b} \in A^+$ such that $\vb{a} \lea \vb{b}$.
    Since $(c,d)\not\in\PhiX (\cb{C}, \vb{a})$, Lemma~\ref{lem:embeddingPhi}
 yields $(c,d)\not\in\PhiX (\cb{C}, \vb{b})$
 and thus $\vb{b}\in\lambda (\cb{C}, (c,d))$. \qed

\section{Relations}    \label{sec:rels}

 A finitary \emph{relation} $R$ on a set $A$ is a subset of $A^I$
 for some finite set $I$. We say a function $f: A^k \rightarrow A$
 \emph{preserves} $R$ if
 $R$ is a subuniverse of $\algop{A}{f}^I$.

 For a clone $\cb{C}$ on a set $A$ and for $m\in\N$, the set of $m$-ary
 functions $\cb{C}^{[m]}$ is a subset of $A^{A^m}$.
 In this sense, a function $f: A^k \rightarrow A$ preserves
 the relation $\cb{C}^{[m]}$ if for all $g_1,\ldots,g_k\in\cb{C}^{[m]}$ the function
\[ A^{m}\rightarrow A,\ x\mapsto f(g_1(x),\ldots,g_k(x)), \] 
 is in $\cb{C}^{[m]}$ again.

 For $\vb{a}\in A^+$ let $|\vb{a}|$ denote the length of $\vb{a}$.

 In the next result we give finitely many relations that determine a clone
 with edge operation.

\begin{thm} \label{thm:relations}
 Let $A$ be a finite set, let $k\in\N$, $k > 1$, let $\cb{C}$ be a clone on $A$
 that contains a $k$-edge operation $t$, and let $\ab{A} := \algop{A}{\cb{C}}$.
 Then the set
     $\{ |\vb{a}| \,\,\setsuchthat \,\, \text{there exists}$
     $(c,d)\in A^2$
    such that 
    $\vb{a}$
     is minimal with respect to $\leq_E$
       in $\lambda (\cb{C}, (c,d)) \}$
 has a supremum $m$ in $\N$, and $\cb{C}$ is the clone of functions that preserve
 the relation $\cb{C}^{[m]}$ and every subuniverse of $\ab{A}^{k-1}$.
\end{thm}

 So by Theorem~\ref{thm:relations} the clone $\cb{C}$ is determined by the
 finitely many relations of arity $\max(|A|^m,k-1)$.
 Apart from the condition on the $m$-ary functions our result resembles the Baker-Pixley
 Theorem (see Theorem 2.1 (5) in~\cite{BP:PIAT}) for clones with near-unanimity operations.

\emph{Proof of Theorem~\ref{thm:relations}:}
 Let $(c,d)\in A^2$.
 Since $(A^+,\leq_E)$ has no infinite antichain by Lemma~\ref{lem:leawpo},
 $\lambda (\cb{C}, (c,d))$ contains only finitely many minimal elements.
 Consequently, as the supremum of finitely many natural numbers, $m$ is finite.
 We note that the set $\{ |\vb{a}| \,\,\setsuchthat \,\, \text{there exists}$
     $(c,d)\in A^2$
    such that 
    $\vb{a}$
     is minimal with respect to $\leq_E$
       in $\lambda (\cb{C}, (c,d)) \}$
 is empty if $\lambda (\cb{C}, (c,d))$ is empty for all $(c,d)\in A^2$. In that case we have
 $m = 1$ as the supremum.

 Let $\cb{D}$ be the clone of functions that preserve $\cb{C}^{[m]}$ and
 every subuniverse of $\ab{A}^{k-1}$.
 Then $\cb{C}\subseteq \cb{D}$ and $\cb{C}^{[m]} = \cb{D}^{[m]}$. We claim that 
\begin{equation} \label{eq:pC}
 \lambda (\cb{C}, (c,d)) \subseteq \lambda (\cb{D}, (c,d)).
\end{equation}
 If $\lambda (\cb{C}, (c,d)) = \emptyset$, the assertion is clear.
 So let $\vb{a}$ be minimal in $\lambda (\cb{C}, (c,d))$.
 Then $(c,d)\not\in\PhiX(\cb{C},\vb{a})$.
 By definition, $m$ is at least the length $|\vb{a}|$ of $\vb{a}$.
 Hence $\cb{C}^{[|\vb{a}|]} = \cb{D}^{[|\vb{a}|]}$, which implies that
 $\PhiX(\cb{C},\vb{a}) = \PhiX(\cb{D},\vb{a})$.
 Thus $\vb{a}\in\lambda (\cb{D}, (c,d))$.
 So we have just proved that every minimal element of $\lambda (\cb{C}, (c,d))$
 is contained in $\lambda (\cb{D}, (c,d))$. Since $\lambda (\cb{C}, (c,d))$ and
 $\lambda (\cb{D}, (c,d))$ are upward closed subsets of the well partially
 ordered set $(A^+,\leq_E)$ by Lemma~\ref{lem:upwardclosed},
 this proves~\eqref{eq:pC}.

 Next we will show that $\cb{D}^{[n]}\subseteq \cb{C}^{[n]}$ for all $n\in\N$.
 For fixed $n\in\N$ and $\vb{a}\in A^n$ we have
\begin{equation} \label{eq:fc}
 \PhiX(\cb{D},\vb{a}) \subseteq\PhiX(\cb{C},\vb{a})
\end{equation}
 by~\eqref{eq:pC}. 

 Note that $F := \cb{C}^{[n]}$ and $G := \cb{D}^{[n]}$ form subuniverses of
 $\ab{A}^{|A|^n}$ with $F\subseteq G$.
 For every $T\subseteq A^n$ with $|T| < k$ we claim that
\begin{equation} \label{eq:FT}
 \pi_T(F) = \pi_T(G).
\end{equation}
 Clearly $\pi_T(F) \subseteq \pi_T(G)$.
 For proving the converse inclusion let $g\in G$, let $l := |T|$, and
 let $T = \{t_1,\ldots, t_l\} = \{ (a_{11},\ldots,a_{1n}), \ldots,
             (a_{l1},\ldots,a_{ln}) \}$.
 We know that $g$ preserves the subuniverse $B$ of
 $\ab{A}^l$ that is generated by $\{ (a_{11},\ldots,a_{l1}),\ldots, 
                                     (a_{1n},\ldots,a_{ln}) \}$.
 From $(g(t_1),\ldots, g(t_l)) \in B$, we obtain
 an $n$-ary term function
 $f$ of $\ab{A}$ such that $ (g (t_1),\ldots, g (t_l)) =
 (f (t_1), \ldots, f(t_l))$. Hence $f|_T = g|_T$, and
 thus $\pi_T(f) = \pi_T(g)$. 
 Hence $\pi_T(F) \supseteq \pi_T(G)$ and we have~\eqref{eq:FT}.
 By~\eqref{eq:fc} and~\eqref{eq:FT} the assumptions of Lemma~\ref{lem:rep}
 are satisfied. Thus $F = G$. 
\qed

     For a finite set $A$ and a set $S$ of finitary relations on $A$,
     we will write $\Pol (A, S)$ for the set of those functions
     on $A$ that preserve all relations in $S$ (cf. \cite{PK:FUR}).

\begin{thm} \label{thm:countable}
        Let $A$ be a finite set, let $k\in\N$, $k>1$, and let 
        $\mathcal{M}_k$ be the set of all clones on $A$ that contain a $k$-edge
 operation.
        Then we have:
        \begin{enumerate}
            \item \label{it:m2} For every 
                                clone $\cb{C}$ in $\mathcal{M}_k$, 
                                there is a finitary
                  relation $R$ on $A$ such that $\cb{C} = \Pol (A, \{ R \})$.
            \item \label{it:m1} There is no infinite descending chain
                  in $(\mathcal{M}_k, \subseteq)$.
            \item \label{it:m3} The set $\mathcal{M}_k$ is finite or 
                  countably infinite.
        \end{enumerate}
    \end{thm}
     \emph{Proof:}
       \eqref{it:m2}  Let $\cb{C}$ be a clone with $k$-edge term on the finite
 set $A$. By Theorem~\ref{thm:relations} there exists a finite set $S$
        of finitary relations on $A$ such that
        $\cb{C} = \Pol (A, S)$. By \cite[p.\ 50]{PK:FUR},
         there is a single finitary relation $R$ on $A$ with
        $\Pol (A, S) = \Pol (A, \{R\})$.

 Now \eqref{it:m1} follows from~\eqref{it:m2} using
        the implication \textsl{(i)'}$\Rightarrow$\textsl{(ii)'} in
        \cite[Charakterisierungssatz 4.1.3]{PK:FUR}.
 
        \eqref{it:m3}
         Every finitary relation on the finite set $A$ is a finite
         subset of the countable set $A^+$. Hence the claim
         follows from~\eqref{it:m2}.
 \qed

\begin{cor} \label{cor:countable}
 Let $A$ be a finite set. Modulo term equivalence, the number of algebras on $A$ that have
 few subpowers is at most countably infinite.
\end{cor}

\emph{Proof:}
 By~\cite[Corollary 3.11]{BI:VWFS} every algebra on $A$ with few subpowers has an edge
 operation in its clone of term functions. 
 Since the number of clones with edge operation on $A$ is at most countably infinite by
 Theorem~\ref{thm:countable}~\eqref{it:m3}, the assertion follows.
\qed

 We recall that a \emph{primitive-positive formula over a language} $\mathcal{R}$ of relation
 symbols is a first-order formula $\varphi(x_1,\dots,x_n)$ of the form
\[ \exists y_1,\dots,y_k\colon (\alpha_1 \wedge \ldots\wedge \alpha_l) \]
 where $\alpha_1,\ldots, \alpha_l$ are atomic formulas, that is, either of the form
 $R(v_1,\dots,v_m)$ for some $R\in\mathcal{R}$ and variables $v_1,\dots,v_m$
 or some equality $v_1 = v_2$ for variables $v_1,v_2$.
 The variables in $\alpha_1,\ldots, \alpha_l$ are from
 $\{ x_1,\dots,x_n \} \cup \{ y_1,\dots,y_k\}$.

 For a set $A$ and $m,n\in\N$, let $R$ be a subset of $A^m$ and let $S$ be a subset of $A^n$.
 We say that $S$ is \emph{primitive-positive definable over} $R$ if there exists a
 primitive-positive formula $\varphi(x_1,\dots,x_n)$ over the language of the relational structure
 $(A,\{R\})$ such that 
\[ (a_1,\dots,a_n)\in S \text{ if and only if } (A,\{R\}) \text{ satisfies } \varphi(a_1,\dots,a_n). \]
 We can now formulate a consequence of Theorem~\ref{thm:countable} that was not known even
 for finite groups $\ab{A}$.

\begin{cor} \label{cor:ppdefinable}
 Let $\ab{A}$ be a finite algebra with few subpowers. Then there exists a subalgebra $R$ of some
 finitary power of $\ab{A}$ such that for every $n\in\N$, every subalgebra $S$ of $\ab{A}^n$ is
 primitive-positive definable over $R$.
\end{cor}

\emph{Proof:}
 By~\cite[Corollary 3.11]{BI:VWFS} the clone $C$ of term operations of $\A$ contains an
 edge operation. So, by Theorem~\ref{thm:countable}~\eqref{it:m2}, we have a finitary relation $R$
 on $A$ such that $C = \Pol(A,\{R\})$.
 Hence by~\cite[Folgerung 1.2.4, Hauptsatz~2.1.3]{PK:FUR}
 every finitary relation $S$ on $A$ that is preserved by all functions in
 $C$ is primitive-positive definable over $R$. Since the finitary relations that are preserved by
 all term functions are exactly the subalgebras of finite powers of $\A$, the result is proved.
\qed

 For the case of finite groups we restate the previous corollary and give some
 explicit bounds on the length of the primitive-positive formula necessary
 to describe an arbitrary relation.

\begin{cor} \label{cor:group}
 Let $\ab{G}$ be a finite, non-trivial group. Then there exists $k\in\N$ and
 a subgroup $H$ of $\ab{G}^k$ with the following property:

    For each $n \in \N$ there are  $l,m \in \N$ with
    $l \le |G|^{n \cdot  \log_2 (|G|)}$ and $m \le l \cdot \log_2 (|G|)$,
    and there is a
    mapping $\sigma : \ul{m} \times \ul{k} \to \ul{l}$
    such that for every subgroup $S$ of $G^n$ there is
    a mapping $\tau : \ul{n} \to \ul{l}$ such that
\[ \begin{array}{ll}
      S = \{  (g_1,\ldots, g_n) \in G^n \setsuchthat
                   \exists a_1,\ldots,a_l \in G : &
                    \bigl( \bigwedge_{i \in \ul{m}}
                         (a_{\sigma(i,1)},\ldots,a_{\sigma(i,k)}) \in H
                    \bigr)  \wedge \\ 
 &                g_1 = a_{\tau(1)} \wedge \ldots \wedge g_n =   
                     a_{\tau(n)} \}.
\end{array} \]
\end{cor}

\emph{Proof:}
 As a subgroup of $\ab{G}^n$, $S$ has a set of generators $\{s_1,\dots,s_e\}$ with
 $e\leq \log_2(|G|^n)$. Let $C$ be the clone of term operations on $\ab{G}$. Then
\begin{equation} \label{eq:S} 
 S = \{ f(s_1,\dots,s_e) \setsuchthat f\in C^{[e]} \}. 
\end{equation}
 By Theorem~\ref{thm:countable}~\eqref{it:m2}, we have $k\in\N$ and some
 subgroup $H$ of $\ab{G}^k$ such that $C$ consists exactly of those functions
 that preserve $H$. In particular
\begin{eqnarray*}
 C^{[e]} & = & \{ f\in G^{G^e} \setsuchthat \bigwedge_{(r_1,\dots,r_e)\in H^e} f(r_1,\dots,r_e) \in H \}, \\
  & = & \bigcap_{(r_1,\dots,r_e)\in H^e} \{ f\in G^{G^e} \setsuchthat f(r_1,\dots,r_e) \in H \}.
\end{eqnarray*}
 Each of the $|H|^e$ many sets in this intersection forms a subgroup of
 $\ab{G}^{G^e}$. So we can choose $\log_2(|G|^{|G|^e})$ many of them whose
 intersection is again equal to $C^{[e]}$. Hence we have $M\subseteq H^e$ with
 $|M|\leq |G|^e \cdot \log_2(|G|)$ such that
\begin{equation} \label{eq:Ce}
 C^{[e]} = \{ f\in G^{G^e} \setsuchthat \bigwedge_{(r_1,\dots,r_e)\in M} f(r_1,\dots,r_e) \in H \}. 
\end{equation}
 Combining~\eqref{eq:S} and~\eqref{eq:Ce} yields
\begin{equation} \label{eq:S2}
 S = \{ g\in G^n \setsuchthat \exists f\in G^{G^e}: \bigwedge_{(r_1,\dots,r_e)\in M} f(r_1,\dots,r_e) \in H \wedge f(s_1,\dots,s_e) = g \}. 
\end{equation} 
 It only remains to rewrite~\eqref{eq:S2}.
 Let $l := |G|^e$, and let $\lambda\colon G^e \to \ul{l}$ be a bijection.
 For $i\in\ul{l}$ define $a_i := f(\lambda^{-1}(i))$.
 Let $m := |M|$, let $\mu\colon \ul{m} \to M$ be a bijection, and let  
 $\sigma \colon \ul{m} \times \ul{k} \to \ul{l}, (i,j) \mapsto \lambda( (\mu(i))_{1j},\dots,(\mu(i))_{ej} )$.
 Note that $l,m$ and $\sigma$ only depend on $n$ but not on $S$.
 Finally define $\tau\colon \ul{n} \to \ul{l}$ by
 $\tau(i) := \lambda(s_{1i},\dots,s_{ei})$. Then the result follows
 from~\eqref{eq:S2}.
\qed

\section{Concluding remarks}

 Using \cite{Id:CCMO} and \cite[Corollary~4.10]{KS:CAWP} together with
 Theorem~\ref{thm:countable}~\eqref{it:m3}, 
 we obtain that the number of clones with $k$-edge term for a fixed integer
 $k > 1$ on a finite set $A$ is finite if $|A| \le 3$, and
 countably infinite if $|A| \ge 4$.

 Given a set $F$ of functions on a finite set $A$ such that $F$ generates a
 clone $C$ with edge operation, Theorem~\ref{thm:countable} guarantees the
 existence of a single relation $R$ that determines $C$;
 however, even if $F$ is finite, it is not yet clear how to find $R$ algorithmically.

 In~\cite{Ko:FSFE} M.~Kozik considered the question whether a function can be 
 obtained as composition of some fixed functions. More precisely, for a
 fixed set of functions $F$ on a finite set $A$ the problem ISTERMFUNCTION
 is the following:

\begin{tabular}{ll}
 INPUT & a function $f: A^n \to A$ \\
 PROBLEM & decide if $f$ is in the clone $C$ on $A$ that is generated by $F$.
\end{tabular}

 He showed that in general this decision problem is EXPTIME-complete.
 If we assume that $F$ contains an edge operation, then there exists 
 some $k$-ary relation $R$ on $A$ such that $C = \Pol (A, \{ R \})$.
 Whether $f$ preserves $R$ can be checked by evaluating $f$ in $k \cdot |R|^n$
 places and performing $|R|^n$ tests whether a given $k$-tuple is an element of $R$.
 Consequently ISTERMFUNCTION is solvable in polynomial time if the algebra
 $\algop{A}{F}$ has few subpowers.

\section{Acknowledgments}
 
 The authors thank J.~Farley and C.~Pech for helpful discussions.
 The second author acknowledges support from Portuguese Project ISFL-1-143 of
 CAUL financed by FCT and FEDER.

\newcommand{\etalchar}[1]{$^{#1}$}
\def\cprime{$'$} \def\cprime{$'$}
\providecommand{\bysame}{\leavevmode\hbox to3em{\hrulefill}\thinspace}
\providecommand{\MR}{\relax\ifhmode\unskip\space\fi MR }
\providecommand{\MRhref}[2]{%
  \href{http://www.ams.org/mathscinet-getitem?mr=#1}{#2}
}
\providecommand{\href}[2]{#2}

\end{document}